# INVERSION OF AN INTEGRAL TRANSFORM AND LADDER REPRESENTATIONS OF $U(1,q)$

John D. Lorch and Lisa A. Mantini

ABSTRACT. An integral transform for $G = U(1,q)$ is studied. The transform maps the positive spin ladder representations of $G$ on a Bargmann-Segal-Fock space $\mathcal{F}_n^{1,q}$ into a space of polynomial-valued functions on the bounded realization $\mathbf{B}^q$ of $G/K$. An inversion is given for the transform and unitary structures are given for the geometric realization of the positive spin ladder representations over $G/K$.

## 0. Introduction

The study of representation theory of linear reductive Lie groups is intimately related to the study of differential equations and the physics that these equations model. Such are the beginnings of the theory concerning the ladder representations of $U(p,q)$ given in this paper.

In [7], Jakobsen and Vergne show the unitarity of a cleverly constructed multiplier representation $\pi_s$ of $U(2,2)$, the 4-1 cover of the conformal group, on a Hilbert space $\mathcal{H}_s$ of holomorphic functions on $U(2,2)/K$ realized as the generalized upper half plane. The elements of $\mathcal{H}_s$ satisfy the mass 0, spin $s$ equations, and the restriction of $\pi_s$ to the Poincaré group is the representation of mass 0, spin $s$.

A previous result that is recurrent throughout [7] is a theorem of Kunze (see [9]) that guarantees the unitarity of certain multiplier representations. The key is that these representation spaces of functions possess a positive definite reproducing kernel that behaves well with respect to the action of the maximal compact subgroup on the image space of these functions.

Meanwhile, the oscillator representation of $U(p,q)$ has been studied and decomposed into irreducible ladder representations in [2], [3], [8], and [16] to name a few. Consider the realization in [3]. Here the oscillator representation $\sigma$ is constructed via the Stone-von Neumann theorem on a Fock space $\mathcal{F}$ of $(p,q)$-holomorphic functions. We use this realization of $\sigma$ in this paper, despite the geometric advantages of the realization given in [2].

1991 *Mathematics Subject Classification.* Primary 22E45, 22E70; Secondary 32L25, 32M15, 58G05, 81R05, 81R25.

*Key words and phrases.* Ladder representations, unitary structures, massless field equations, twistor theory, Penrose transform.

The first author gratefully acknowledges Dissertation Fellowship support from the Mathematics Department of Oklahoma State University.

The second author gratefully acknowledges post-doctoral fellowship support from the American Association of University Women, sabbatical year support from Oklahoma State University, and support from NSF Grant DMS 9304580 at the Institute for Advanced Study.

This paper is in final form and no version of it will be submitted for publication elsewhere.

Typeset by $\mathcal{A}_{\mathcal{M}}\mathcal{S}$-TEX

The relationship between the representations $\pi_s$ and the oscillator representation of $U(2,2)$ noted in [7] is generalized by Mantini in [10]. It is in [10] that we see the positive spin ladder representations $\mathcal{F}_n$ of $U(p,q)$ realized via an integral transform $\Phi_n$ as a subspace $\mathcal{S}_n$ of sections of a homogeneous holomorphic vector bundle $\mathbf{E}_n$ over $U(p,q)/K$ in its model as a generalized unit disk $\mathbf{D}$. The fiber of $\mathbf{E}_n$ over any point $\zeta \in \mathbf{D}$ is essentially the space $\bar{\mathcal{P}}(n, \mathbf{C}^q)$ of antiholomorphic polynomials homogeneous of degree $n$ in $q$ complex variables. If the real rank of $U(p,q)$ is larger than one, or if the rank of $\mathbf{E}_n$ is larger than one, these sections satisfy certain linear partial differential equations. In the $U(2,2)$ case, these equations are the massless field equations and the natural action on $\mathcal{S}_n$ is equivalent to the multiplier representation of [7]. The transform $\Phi_n$ bears many similarities to the Penrose transform (see [12] and [17]). We will discuss a version of $\Phi_n$ more carefully in section 2 below. There also exists an analogous construction of the negative spin ladder representations (see [11]), but here we restrict to the positive spin case.

Since $\Phi_n$ is one-to-one on $\mathcal{F}_n$ for $n \geq 0$, $\Phi_n$ is invertible. In section 3, an inverse for $\Phi_n$ in the case when $G = U(1,q)$ is presented. Trailing along behind the inversion of $\Phi_n$ is the machinery that allows the production of a Hilbert space structure on $\mathcal{S}_n$ so that $\mathcal{S}_n$ is unitarily equivalent to $\mathcal{F}_n$ via $\Phi_n$. Such structures are exhibited in section 4 for this case.

Motivation for our results is essentially two-fold. On one hand, we consider the results of Davidson [3] (which again depend critically upon results of Kunze [9]). He gives natural unitary norms on $U(p,q)/K$ for discrete series representations of $U(p,q)$. These structures take on the form

$$(f_1, f_2) = \int_{\mathbf{D}} \langle Q_\lambda(\zeta,\zeta)^{-1} f_1(\zeta), f_2(\zeta)\rangle \, d\mu(\zeta),$$

where $f_i \in \mathcal{H}(\mathbf{D}, \lambda)$ are polynomial-valued functions on the disk, $d\mu$ is the $U(p,q)$-invariant measure on $\mathbf{D}$, and $Q_\lambda$ is the reproducing kernel for $\mathcal{H}(\mathbf{D}, \lambda)$. The inner product $\langle \cdot, \cdot \rangle$ consists of an integral with respect to the appropriate Gaussian measure. For the ladder representations of $U(p,q)$ which are not square-integrable, an analogous norm does not converge (see Lemma 4.2 of [3]).

On the other hand, we are also inspired by inversion formulas for the Radon transform, as given in [5]. Let $\mathbf{P}^n$ denote the space of all hyperplanes in $\mathbf{R}^n$, and let $f$ be a function on $\mathbf{R}^n$ that is integrable on each element $\xi$ of $\mathbf{P}^n$ with respect to Euclidean measure $dm_\xi$ on $\xi$. The Radon transform $\hat{f}$ of $f$ is then given by

$$\hat{f}(\xi) = \int_\xi f(x) \, dm_\xi(x),$$

while the dual transform for continuous functions $\varphi$ on $\mathbf{P}^n$ is given by

$$\check{\varphi}(x) = \int_{\xi_x} \varphi(\xi) \, d\mu_x(\xi),$$

where $\xi_x = \{\xi \in \mathbf{P}^n \mid x \in \xi\}$ and $d\mu_x$ is the measure of total mass one on $\xi_x$ that is invariant under rotation about $x$. A Schwartz function $f$ can be recovered from its Radon transform via the formula

$$cf = \mathcal{L}^{\frac{n-1}{2}}((\hat{f})\check{\;}),$$

where $c$ is a constant and $\mathcal{L}$ is the Laplacian on $\mathbf{R}^n$. If $n$ is even, $\mathcal{L}^{\frac{n-1}{2}}$ is defined in terms of Riesz potentials (see Lemma 2.35 of [5]). Operators playing a role similar to these fractional powers of the Laplacian appear in our results.

In large part, we achieve our results by constructing two special operators on $\mathcal{O}(\mathbf{D}, \bar{\mathcal{P}}(n, \mathbf{C}^q))$, $P$ and $L$. Let $\phi$ be the image of a monomial in $\mathcal{F}_n$ under $\Phi_n$. Fixing $\zeta \in \mathbf{D}$ and considering $\phi(\zeta)$ as an element of $\bar{\mathcal{P}}(n, \mathbf{C}^q)$, the mapping $P$ is a projection operator such that $P\phi(\zeta)$ is the nonzero summand of $\phi(\zeta)$ of highest possible weight with respect to the $K$ action on $\bar{\mathcal{P}}(n, \mathbf{C}^q)$. The mapping $L$ is simply a weighted version of the mapping $P$ which enables the overall inversion and inner products to become unitary. Using this mapping $L$, we produce an inversion formula for $\Phi_n$ and the desired unitary structures. These have the form

$$(\Phi_n^{-1}\phi)(z) = \lim_{t \to 1^-} \frac{d^q}{dt^q}\left[t^q \int_{\mathbf{D}} \int_{\mathbf{C}^q} (L\phi)(t\zeta, v) e^{(z_1\bar{\zeta} + v^T)\bar{z}_S} \, e^{-|v|^2} \, dm(v) dm(\zeta)\right],$$

and

$$((\phi_1, \phi_2)) = \lim_{t \to 1^-} \frac{d^q}{dt^q}\left[t^q \int_{\mathbf{D}} \langle (L\phi_1)(t\zeta, \cdot), \phi_2(\zeta, \cdot) \rangle \, dm(\zeta)\right],$$

respectively.

Alternative unitary structures for most of the ladder representations of $U(p, q)$ are given by Rawnsley, Schmid, and Wolf (see section 13 of [14]). These representations occur in $L^2$-cohomology of certain homogeneous, holomorphic vector bundles over $U(p, q)/(U(1) \times U(p-1, q))$.

The authors would like to thank M. Davidson, E. Dunne, E. Grinberg, and R. Zierau for helpful conversations. We are grateful to the referee for some useful remarks. Finally, we would like to thank the editors of this volume for their efforts on behalf of the authors as well as for their patience and support.

## 1. Preliminaries

In this section we wish to gain familiarity with the oscillator representation of $U(p, q)$ realized on a certain Fock space. The generalized unit disk and an associated representation will be defined. We first introduce the multi-index notation that will be prevalent throughout the paper.

### 1.1. Notation.

(a) For $r \in \mathbf{N}$, let $\mathbf{N}_0$ denote the set of nonnegative integers and let $\mathbf{N}_0^r$ denote the set of $r$-tuples of nonnegative integers.

(b) If $m = (m_1, m_2, \ldots, m_r) \in \mathbf{N}_0^r$, then $m! := m_1! m_2! \cdots m_r!$ and $|m| := m_1 + m_2 + \cdots + m_r$. If $z \in \mathbf{C}^r$, then $z^m := z_1^{m_1} z_2^{m_2} \cdots z_r^{m_r}$.

(c) Let $\alpha \in \mathbf{N}_0^s$, $1 \leq s \leq r$, and let $J = (j_1, j_2, \ldots, j_s) \in \mathbf{N}_0^s$ be such that $1 \leq j_1 < j_2 < \cdots < j_s \leq r$. Then, if $z \in \mathbf{C}^r$, $z_J := (z_{j_1}, z_{j_2}, \ldots, z_{j_s})$, $|z_J|^2 := |z_{j_1}|^2 + |z_{j_2}|^2 + \cdots + |z_{j_s}|^2$, and $z_J^\alpha := z_{j_1}^{\alpha_1} z_{j_2}^{\alpha_2} \cdots z_{j_s}^{\alpha_s}$.

(d) Let $R_p := (1, 2, \ldots, p)$ and $S_{p,q} := (p+1, p+2, \ldots, p+q)$. When there is no danger of confusion, we will let $R$ and $S$ represent $R_p$ and $S_{p,q}$, respectively.

(e) For $n \in \mathbf{N}_0$ and $q \in \mathbf{N}$, let $\mathbf{N}_0^q(n) := \{\eta \in \mathbf{N}_0^q \mid |\eta| = n\}$, and let $d(n, q)$ denote the cardinality of $\mathbf{N}_0^q(n)$. Note that $d(n, q) = \dim \bar{\mathcal{P}}(n, \mathbf{C}^q)$.

(f) The symbol $\prec$ will denote lexicographic, or dictionary, order on $\mathbf{N}_0^q(n)$. Here $\alpha \prec \beta$ if there exists $j \in \{1, \ldots, q\}$ such that $\alpha_i = \beta_i$ for $i < j$ and $\alpha_j < \beta_j$. In this case we say that $\beta$ succeeds $\alpha$.

(g) We have a partial order $\leq$ on $\mathbf{N}_0^q$ given by $\alpha \leq \beta$ if $\alpha_i \leq \beta_i$ for $i \in \{1, \ldots, q\}$.

(h) Let $\epsilon_j^q := (0, \ldots, 0, 1, 0, \ldots, 0) \in \mathbf{N}_0^q$, where the support lies in the $j$-th component of $\epsilon_j^q$. We will write $\epsilon_j$ for $\epsilon_j^q$ when there is no danger of confusion.

(i) For $p, q \in \mathbf{N}$, we let $M^{p \times q}(\mathbf{C})$ denote the set of $p \times q$ matrices with complex entries. If $X \in M^{p \times q}(\mathbf{C})$, the conjugate transpose of $X$ is written as $X^*$.

Recall that for $p, q \in \mathbf{N}$, $U(p, q) = \{g \in GL(p+q, \mathbf{C}) \mid g^* I_{p,q} g = I_{p,q}\}$, where $I_{p,q}$ is the element of $M^{(p+q) \times (p+q)}(\mathbf{C})$ given by

$$I_{p,q} = \begin{pmatrix} I_p & 0 \\ 0 & -I_q \end{pmatrix}.$$

If $g \in G = U(p, q)$, the matrix $g$ is occasionally written in block form as

$$g = \begin{pmatrix} A & B \\ C & D \end{pmatrix},$$

where $A$ and $D$ are $p \times p$ and $q \times q$ matrices, respectively. In addition, the subgroup $K \simeq U(p) \times U(q)$ of $U(p, q)$ given by

$$K := \left\{ \begin{pmatrix} A & 0 \\ 0 & D \end{pmatrix} \mid A \in U(p), \ D \in U(q) \right\}$$

is a maximal compact subgroup.

Let $dm(z)$ denote the normalization of Lebesgue measure on $\mathbf{C}^r$ given by

$$dm(z) := \frac{1}{\pi^r} dx_1 dy_1 \ldots dx_r dy_r.$$

The following lemmas will be useful.

**1.2. Lemma.** *Let $k, l \in \mathbf{N}_0$, and let $\delta_{kl} = 1$ if $k = l$ and $\delta_{kl} = 0$ otherwise. Then*

$$\int_{\mathbf{C}} z^k \bar{z}^l e^{-|z|^2} \, dm(z) = \delta_{kl} \, k!.$$

**1.3. Lemma.** *Suppose that $U \in \mathbf{C}^{q \times q}$ is such that $U + U^* \gg 0$, and $\phi$ is a holomorphic function of $z \in \mathbf{C}^q$ for which $\int \phi(z) \exp(-z^* U z) \, dm(z)$ converges absolutely. Then*

$$\int_{\mathbf{C}^q} \phi(z) e^{-z^* U z} \, dm(z) = \frac{\phi(0)}{\det U}.$$

**1.4. Lemma.** *Suppose that $q \in \mathbf{N}$ and $w \in \mathbf{C}$. The series*

$$\sum_{m \in \mathbf{N}_0^q} \frac{w^{|m|}}{|m|^q}$$

*converges if $|w| < 1$ and diverges if $|w| > 1$.*

Finally, we say that a complex-valued function $f$ on $\mathbf{C}^{p+q}$ is *(p,q)-holomorphic* if $f$ is holomorphic in $z_R$ and antiholomorphic in $z_S$.

**1.1. The Oscillator Representation of** $U(p,q)$**.** Consider the Bargmann-Segal-Fock space

$$\mathcal{F}^{p,q} = \{f : \mathbf{C}^{p+q} \to \mathbf{C} \mid f \text{ is } (p,q)\text{-holomorphic}$$
(1.1)
$$\text{and} \int_{\mathbf{C}^{p+q}} |f(z)|^2 e^{-|z|^2} dm(z) < \infty\}.$$

By polarizing the integral in (1.1), one obtains the inner product for the Hilbert space $\mathcal{F}^{p,q}$, which will be denoted by $\langle \cdot, \cdot \rangle$. Let $f \in \mathcal{F}^{p,q}$ and suppose

$$f(z) = \sum_{\substack{l \in \mathbf{N}_0^p \\ m \in \mathbf{N}_0^q}} a_{l,m} \, z_R^l \, \bar{z}_S^m.$$

It follows from Lemma 1.2 that

(1.2)
$$\langle f, f \rangle = \sum_{\substack{l \in \mathbf{N}_0^p \\ m \in \mathbf{N}_0^q}} |a_{l,m}|^2 \, l! \, m!.$$

As an application of Lemma 1.3, we obtain the reproducing formula for $f \in \mathcal{F}^{p,q}$ given by

(1.3)
$$f(z) = \int_{\mathbf{C}^{p+q}} f(w) e^{w_R^* z_R + z_S^* w_S} e^{-|w|^2} \, dm(w).$$

We now define a version of the oscillator representation $\sigma$ of $G$.

**1.5. Definition.** For $g = \begin{pmatrix} A & B \\ C & D \end{pmatrix} \in G$ and $f \in \mathcal{F}^{p,q}$, define $\sigma(g)f$ by

$$[\sigma(g)f](z) = \det(D^*) \int_{\mathbf{C}^{p+q}} f(g^{-1}w) e^{-\frac{1}{2}|g^{-1}w|^2} e^{w_R^* z_R + z_S^* w_S} e^{-\frac{1}{2}|w|^2} dm(w).$$

It is not immediate that $\sigma$ defines a unitary representation of $G$ on $\mathcal{F}^{p,q}$. Nonetheless we have the following theorem, due to Blattner and Rawnsley.

**1.6. Theorem** [2]**.** *The map $\sigma$ of Definition 1.5 is a unitary representation of $G$ on $\mathcal{F}^{p,q}$.*

Now, for $n \in \mathbf{Z}$, let

$$\mathcal{F}_n^{p,q} := \{f \in \mathcal{F}^{p,q} \mid f(e^{-i\theta}z) = e^{in\theta}f(z), \theta \in \mathbf{R}\}.$$

A short computation shows that $\mathcal{F}_n^{p,q}$ is a $G$-invariant subspace of $\mathcal{F}^{p,q}$. For $g \in G$, let $\sigma_n(g)$ denote the restriction of $\sigma(g)$ to $\mathcal{F}_n^{p,q}$.

**1.7. Theorem** [16]**.** *Let $n \in \mathbf{Z}$. The representation $\sigma_n$ of $G$ on $\mathcal{F}_n^{p,q}$ is irreducible, and there is an orthogonal direct sum decomposition*

$$\mathcal{F}^{p,q} = \bigoplus_{n \in \mathbf{Z}} \mathcal{F}_n^{p,q}.$$

*Remark.* The representations $\sigma_n$ are often referred to as *ladder* representations, since the highest weights of their $K$-types lie equally spaced along a ray in the weight lattice.

**1.2. The Generalized Unit Disk.** For $p, q \in \mathbf{N}$, let $\mathbf{C}^{p,q}$ denote $\mathbf{C}^{p+q}$ endowed with the hermitian form $h$ with matrix $I_{p,q}$. A subspace $V$ of $\mathbf{C}^{p,q}$ is *negative* if the restriction of $h$ to $V$ is negative definite. Let $\mathbf{M}_{p,q}^-$ denote the set of all negative $q$-planes in $\mathbf{C}^{p,q}$. The *generalized unit disk* $\mathbf{D}_{p,q}$ is given by

$$(1.4) \qquad \mathbf{D}_{p,q} = \{\zeta \in M^{p\times q}(\mathbf{C}) \mid I_q - \zeta^*\zeta \gg 0\}.$$

**1.8. Proposition.** *The mapping $\zeta \mapsto V(\zeta)$ of $\mathbf{D}_{p,q}$ into $\mathbf{M}_{p,q}^-$ where*

$$V(\zeta) = \text{col span} \begin{bmatrix} \zeta \\ I_q \end{bmatrix}$$

*gives a parameterization of $\mathbf{M}_{p,q}^-$ by $\mathbf{D}_{p,q}$.*

*Proof.* The mapping is clearly well-defined and injective into the space of $q$-planes. The $q$-plane $V(\zeta)$ associated to $\zeta \in \mathbf{D}_{p,q}$ is negative due to the positivity condition in (1.4). It remains to show that the mapping is surjective. Let

$$V = \text{col span} \begin{bmatrix} U_1 \\ U_2 \end{bmatrix} \in \mathbf{M}_{p,q}^-,$$

where $U_1$ and $U_2$ are $p \times q$ and $q \times q$ matrices, respectively. Now $U_2$ must be nonsingular, as otherwise $V$ will admit positive vectors. Hence we can change basis to see that $V = V(U_1 U_2^{-1})$. Since $V$ is negative, $U_1 U_2^{-1}$ satisfies the conditions of (1.4), therefore $U_1 U_2^{-1} \in \mathbf{D}_{p,q}$. ∎

The natural left action of $G$ on $\mathbf{M}_{p,q}^-$ given by

$$g.V(\zeta) = \text{col span}\left(g \begin{bmatrix} \zeta \\ I_q \end{bmatrix}\right)$$

is transitive, and the stabilizer of $V(0)$ is $K$. By Proposition 1.8, $\mathbf{D}_{p,q}$ is diffeomorphic to $G/K$ and there is a natural action of $G$ on $\mathbf{D}_{p,q}$ that is stabilized at $0 \in \mathbf{D}_{p,q}$ by $K$. This action is given by

$$g.\zeta = (A\zeta + B)(C\zeta + D)^{-1},$$

where $g = \begin{pmatrix} A & B \\ C & D \end{pmatrix} \in G$ and $\zeta \in \mathbf{D}_{p,q}$.

Note that $\mathbf{D}_{p,q}$ is a bounded domain in $\mathbf{C}^{pq}$. In fact, it is one of the four *classical domains* (cf. [6] or [13]). In the case we are concerned with, that is $p = 1$, $\mathbf{D}_{1,q}$ is simply the unit ball $\mathbf{B}^q$ in $\mathbf{C}^q$.

We will rely on our ability to integrate over $\mathbf{B}^q$, and toward that end we use the following result. For details, see [15].

**1.9. Lemma.** *Let $\eta, \gamma \in \mathbf{N}_0^q$ and suppose that $|\eta| = n$ and $|\gamma| = m$. Then*

$$\int_{\mathbf{B}^q} \zeta^\eta \bar{\zeta}^\gamma dm(\zeta) = \delta_{\eta,\gamma}\left(\frac{\eta!}{(n+q)!}\right),$$

where $\delta_{\eta,\gamma} = 1$ if $\eta = \gamma$ and $\delta_{\eta,\gamma} = 0$ *otherwise.*

Let $\bar{\mathcal{P}}(n, \mathbf{C}^q)$ denote the set of antiholomorphic polynomials on $\mathbf{C}^q$ that are homogeneous of degree $n$ and let $\mathcal{O}(\mathbf{D}_{p,q}, \bar{\mathcal{P}}(n, \mathbf{C}^q))$ denote the set of functions holomorphic on $\mathbf{D}_{p,q}$ taking values in $\bar{\mathcal{P}}(n, \mathbf{C}^q)$. We conclude this section with the construction of a representation $\omega_n$ of $G$ on $\mathcal{O}(\mathbf{D}_{p,q}, \bar{\mathcal{P}}(n, \mathbf{C}^q))$. Define, for $n \geq 0$, a mapping $J_n : G \times \mathbf{D}_{p,q} \to GL(\bar{\mathcal{P}}(n, \mathbf{C}^q))$ by

$$(1.5) \qquad J_n(g, \zeta)f(v) = \det[C\zeta + D]f([C\zeta + D]^*v)$$

for $g = \begin{pmatrix} A & B \\ C & D \end{pmatrix}$, $\zeta \in \mathbf{D}_{p,q}$ and $v \in \mathbf{C}^q$. It is easily shown that $J_n$ is a *factor of automorphy*. That is, if $g_1, g_2 \in G$ and $\zeta \in \mathbf{D}_{p,q}$ then $J_n(I, \zeta) = \mathrm{Id}$ and $J_n(g_1 g_2, \zeta) = J_n(g_1, g_2.\zeta) \circ J_n(g_2, \zeta)$.

**1.10. Definition.** Suppose $g \in G$ and $\phi \in \mathcal{O}(\mathbf{D}_{p,q}, \bar{\mathcal{P}}(n, \mathbf{C}^q))$. Define $\omega_n(g)\phi$ by

$$(\omega_n(g)\phi)(\zeta, v) = J_n(g^{-1}, \zeta)^{-1}\phi(g^{-1}.\zeta, v)$$

for $\zeta \in \mathbf{D}_{p,q}$ and $v \in \mathbf{C}^q$.

Note that $\omega_n$ is a multiplier representation of the type discussed in [9]. This representation corresponds to the natural geometric action of $G$ on the vector bundle $\mathbf{E}_n$. In the next section we will intertwine $\sigma_n$ and $\omega_n$ via an integral transform.

## 2. The Integral Transform

In [10], Mantini constructs an integral transform $\Phi_n$ for each $n \in \mathbf{Z}$, $n \geq 0$. We may think of $\Phi_n$ as a mapping of $\mathcal{F}_n^{p,q}$ into $\mathcal{O}(\mathbf{D}_{p,q}, \bar{\mathcal{P}}(n, \mathbf{C}^q))$. This transform is analogous to the Penrose transform (see [12]). In case either the real rank of $G$ is larger than one or the dimension of $\bar{\mathcal{P}}(n, \mathbf{C}^q)$ is larger than one (i.e., when $n > 0$ and $q > 1$), the elements in the image of $\Phi_n$ satisfy certain linear partial differential equations. These equations generalize the massless field equations of mathematical physics which hold when $G = U(2,2)$ (see [10] or [7]). The image of the map $\Phi_n$ in [10] is in a space of sections of a holomorphic, homogeneous vector bundle on which the group acts by the natural, geometric action. Thus $\Phi_n$ implicitly gives a geometric construction of the ladder representations $\sigma_n$, although the inner product structure on this realization has previously been unknown.

At this point let us recall that our primary goal is to create an inversion for $\Phi_n$ and then use this inversion to develop unitary structures for the geometric realization of $\sigma_n$ obtained via $\Phi_n$. In this section we familiarize ourselves with $\Phi_n$ and some of its important properties. Further details may be found in [10] and [11].

**2.1. Proposition [10].** *Suppose $n \geq 0$ and $f \in \mathcal{F}_n^{p,q}$. Then we have an integral transform $\Phi_n : \mathcal{F}_n^{p,q} \to \mathcal{O}(\mathbf{D}_{p,q}, \bar{\mathcal{P}}(n, \mathbf{C}^q))$ given by*

$$(\Phi_n f)(\zeta, v) = \int_{\mathbf{C}^q} f(\zeta w, w) e^{v^* w} e^{-|w|^2} \, dm(w),$$

*where $v \in \mathbf{C}^q$ and $(\Phi_n f)(\zeta, v)$ depends holomorphically on $\zeta \in \mathbf{D}_{p,q}$.*

The transform of Proposition 2.1 is simplified from the transform constructed by Mantini in [10]. Here we have used the trivialization of the vector bundle in [10] with the multiplier action $\omega_n$. We have also used Davidson's model of the Fock space representation $\sigma_n$ instead of the $L^2$-cohomology realization of Blattner and Rawnsley (see [2]). The transform is essentially given by restriction to a negative $q$-plane determined by $\zeta \in \mathbf{D}_{p,q}$ followed by a projection operator, giving rise to an $L^2$ version of the Penrose transform. Notice that the projection operator is the specialization to signature $(0,q)$ of the reproducing formula given in (1.3).

The first important feature of $\Phi_n$ is that it serves as an intertwining map between the representations $\sigma_n$ and $\omega_n$.

**2.2. Theorem** [10]. *The mapping $\Phi_n$ is one-to-one for $n \geq 0$ and identically zero for $n < 0$. Furthermore, for all $f \in \mathcal{F}_n^{p,q}$, $g \in G$, and $n \geq 0$, the mapping $\Phi_n$ satisfies*
$$\omega_n(g)(\Phi_n f) = \Phi_n(\sigma_n(g)f).$$

Theorem 2.2 implies that for $n \geq 0$, $\Phi_n$ has an inverse defined on $\Phi_n(\mathcal{F}_n^{p,q})$. Also, it says that $\sigma_n$ is realized on $\Phi_n(\mathcal{F}_n^{p,q})$ by a subrepresentation of $\omega_n$.

**2.3. Lemma** [10]. *Fix $n \in \mathbf{Z}$ and $\zeta \in \mathbf{D}_{p,q}$. The mapping of $\mathcal{F}_n^{p,q}$ into $\mathcal{F}_n^{0,q}$ given by*
$$f \mapsto \Phi_n f(\zeta, \cdot)$$
*is continuous.*

The next thing to note is that for many values of $p$ and $q$, the elements of $\Phi_n(\mathcal{F}_n^{p,q})$ satisfy certain interesting differential equations. Using the notation of [4], let $\square_{p,q}$ be the $(p+1) \times q$ matrix of differential operators where
$$(\square_{p,q})_{ij} = \begin{cases} \frac{\partial}{\partial \zeta_{ij}}, & \text{for } 1 \leq i \leq p \text{ and } 1 \leq j \leq q, \\ \frac{\partial}{\partial \bar{v}_j}, & \text{for } i = p+1 \text{ and } 1 \leq j \leq q. \end{cases}$$

**2.4. Theorem** [10].
  (a) *Let $p, q \in \mathbf{N}$ be such that both $p$ and $q$ are greater than one. Then $\Phi_n(\mathcal{F}_n^{p,q})$ is a subset of $\mathcal{O}(\mathbf{D}_{p,q}, \bar{\mathcal{P}}(n, \mathbf{C}^q))$ annihilated by the differential operators formed by taking $2 \times 2$-minors of $\square_{p,q}$.*
  (b) *If $p = 1$, $n \geq 1$ and $q > 1$, then $\Phi_n(\mathcal{F}_n^{1,q})$ is a subset of $\mathcal{O}(\mathbf{B}^q, \bar{\mathcal{P}}(n, \mathbf{C}^q))$ that is annihilated by the differential operators formed by taking $2 \times 2$-minors of $\square_{1,q}$.*

*Remark.* In case $p = q = 2$, the differential operators formed by taking $2 \times 2$-minors of $\square_{2,2}$ are the complexifications of the massless field equations. Also, it has been shown in [4] that the $K$-finite vectors in $\mathcal{O}(\mathbf{D}_{p,q}, \bar{\mathcal{P}}(n, \mathbf{C}^q))$ with respect to $\omega_n$ that satisfy the differential equations mentioned in the theorem are exactly the $K$-finite vectors in $\Phi_n(\mathcal{F}_n^{p,q})$.

Suppose that $n \in \mathbf{N}_0$ and $\phi \in \mathcal{O}(\mathbf{D}_{p,q}, \bar{\mathcal{P}}(n, \mathbf{C}^q))$. Then there exists, for each $\eta \in \mathbf{N}_0^q(n)$, a holomorphic complex-valued function $\psi^{(\eta)}$ on $\mathbf{D}_{p,q}$ such that

$$(2.1) \qquad \phi(\zeta, v) = \sum_{\eta \in \mathbf{N}_0^q(n)} \psi^{(\eta)}(\zeta) \frac{\bar{v}^\eta}{\eta!}.$$

The component functions $\psi^{(\eta)}$ in (2.1) are indexed by $\eta \in \mathbf{N}_0^q(n)$ under the lexicographic order $\prec$ (see part (f) of Notation 1.1), with indices ranging in order from $(0, \ldots, 0, n)$ through $(n, 0, \ldots, 0)$. In what follows, this description of $\phi(\zeta, v)$ will sometimes be convenient.

We finish this section with a lemma that will be of importance later on. Define, for $\alpha \in \mathbf{N}_0^p$, $\beta \in \mathbf{N}_0^q$ and $z \in \mathbf{C}^{p+q}$,

$$f_{\alpha\beta}(z) := z_R^\alpha \bar{z}_S^\beta. \tag{2.2}$$

We wish to compute the image of $f_{\alpha\beta}$ under the transform $\Phi_{(|\beta|-|\alpha|)}$. In order to simplify our expressions, we first introduce some notation.

Given $\gamma \in M^{p \times q}(\mathbf{C})$ and $i \in \{1, 2, \ldots, p\}$, let $\gamma_{(i)}$ denote the $i$-th row of $\gamma$. Let $r(\gamma)$ denote the element of $\mathbf{C}^q$ obtained by taking the sum of the rows of $\gamma$, and let $c(\gamma) \in \mathbf{C}^p$ denote the sum of the columns of $\gamma$. Given $\alpha \in \mathbf{N}_0^p$ and $\tilde{\alpha} \in \mathbf{N}_0^q$, let $M^{p \times q}(\alpha, \tilde{\alpha})$ denote the subset of $M^{p \times q}(\mathbf{N}_0)$ whose elements $\gamma$ satisfy

$$c(\gamma) = \alpha \quad \text{and} \quad r(\gamma) = \tilde{\alpha}.$$

Notice that $M^{p \times q}(\alpha, \tilde{\alpha})$ is nonempty only when $|\alpha| = |\tilde{\alpha}|$. Now, if $\gamma \in M^{p \times q}(\mathbf{N}_0)$ and $\zeta \in \mathbf{D}_{p,q}$, we may write $\zeta^\gamma$ and $\gamma!$ for

$$\prod_{i=1}^p \zeta_{(i)}^{\gamma_{(i)}} \quad \text{and} \quad \prod_{i=1}^p \gamma_{(i)}!,$$

respectively. Finally, when $p$ and $q$ are understood, we let $M(\alpha, \tilde{\alpha})$ stand for $M^{p \times q}(\alpha, \tilde{\alpha})$.

**2.5. Lemma.** *Suppose that $n \geq 0$, $\alpha \in \mathbf{N}_0^p$ and $\beta \in \mathbf{N}_0^q$ with $|\beta| - |\alpha| = n$. Then $f_{\alpha\beta} \in \mathcal{F}_n^{p,q}$ and*

$$(\Phi_n f_{\alpha\beta})(\zeta, v) = \sum_{\substack{\eta \in \mathbf{N}_0^q(n) \\ \eta \leq \beta}} \sum_{\gamma \in M(\alpha, \beta - \eta)} \left(\frac{\alpha! \, \beta!}{\gamma!} \zeta^\gamma\right) \frac{\bar{v}^\eta}{\eta!},$$

*where $\zeta \in \mathbf{D}_{p,q}$ and $v \in \mathbf{C}^q$. The order $\leq$ is as in part (g) of Notation 1.1.*

*Proof.* Let $\phi$ denote $\Phi_n f_{\alpha\beta}$. According to Proposition 2.1, we have

$$\phi(\zeta, v) = \int_{\mathbf{C}^q} f_{\alpha\beta}(\zeta w, w) e^{v^* w} e^{-|w|^2} \, dm(w)$$

$$= \int_{\mathbf{C}^q} \left[\prod_{i=1}^p \left(\sum_{j=1}^q \zeta_{ij} w_j\right)^{\alpha_i}\right] \bar{w}^\beta e^{v^* w} e^{-|w|^2} \, dm(w).$$

At this point, we apply the generalized binomial theorem and obtain

$$\phi(\zeta, v) = \int_{\mathbf{C}^q} \left[\prod_{i=1}^p \Big(\sum_{\substack{\gamma_{(i)} \in \mathbf{N}_0^q \\ |\gamma_{(i)}| = \alpha_i}} \zeta_{(i)}^{\gamma_{(i)}} w^{\gamma_{(i)}} \frac{\alpha_i!}{\gamma_{(i)}!}\Big)\right] \bar{w}^\beta e^{v^* w} e^{-|w|^2} \, dm(w)$$

$$= \sum_{\substack{\gamma_{(1)}, \ldots, \gamma_{(p)} \in \mathbf{N}_0^q \\ |\gamma_{(i)}| = \alpha_i}} \frac{\alpha! \, \zeta^\gamma}{\gamma!} \int_{\mathbf{C}^q} w^{r(\gamma)} \bar{w}^\beta e^{v^* w} e^{-|w|^2} \, dm(w),$$

where $\gamma$ is the $p \times q$ matrix with rows $\gamma_{(1)}, \ldots, \gamma_{(p)}$.

Now, we are able to compute the integral in the preceding line by expanding the exponential function into its Taylor series and then applying Lemma 1.2 repeatedly. Doing so, we obtain

$$\phi(\zeta, v) = \sum_{\substack{\gamma \in M^{p \times q}(\mathbf{N}_0) \\ c(\gamma) = \alpha}} \frac{\alpha! \, \zeta^{\gamma}}{\gamma!} \, \epsilon(\beta, \gamma) \left[ \frac{\bar{v}^{\beta - r(\gamma)} \beta!}{(\beta - r(\gamma))!} \right],$$

where $\epsilon(\beta, \gamma) = 1$ if $\beta \geq r(\gamma)$, and $\epsilon(\beta, \gamma) = 0$ if $\beta - r(\gamma)$ has any negative components. The result then follows immediately by gathering terms of the form $\bar{v}^{\eta}/\eta!$ for each $\eta \in \mathbf{N}_0^q(n)$. ∎

**2.6. Corollary.** *Fix $p = 1$, $n \geq 0$ and $\nu \in \mathbf{N}_0^q$ with $|\nu| \geq n$. Let $f_{\nu,n} \in \mathcal{F}_n^{1,q}$ be defined by*

$$f_{\nu,n}(z) := z_1^{|\nu| - n} \bar{z}_S^{\nu},$$

*and let $\phi_\nu = \Phi_n(f_{\nu,n})$. Then*

$$\phi_\nu(\zeta, v) = \sum_{\substack{\eta \in \mathbf{N}_0^q(n) \\ \eta \leq \nu}} \left( \frac{\nu! \, (|\nu| - n)!}{(\nu - \eta)!} \zeta^{\nu - \eta} \right) \frac{\bar{v}^{\eta}}{\eta!}.$$

## 3. An Inversion Formula

The goal of this section is to give an inversion formula for the integral transform $\Phi_n$ of section 2, when $n \geq 0$, $p = 1$, and $q \geq 1$. Throughout this section, $p$, $q$, and $n$ will remain as in the previous sentence. We begin with the introduction of two mappings of $\mathcal{O}(\mathbf{B}^q, \bar{\mathcal{P}}(n, \mathbf{C}^q))$ into itself. One mapping, $P$, is a projection that isolates a nonzero summand of $\Phi_n f_{\nu,n}$. The other mapping, $L$, is simply a weighted version of $P$ which we then use to invert the transform $\Phi_n$. The weighting of $L$ is chosen to ensure the $G$-equivariance of our inversion formula.

Let $\nu \in \mathbf{N}_0^q$ with $|\nu| \geq n$ and let $f_{\nu,n}(z) = z_1^{|\nu| - n} \bar{z}_S^{\nu}$, as in Corollary 2.6. For each $\eta \in \mathbf{N}_0^q(n)$ we define a complex-valued function $\psi_\nu^{(\eta)}$ on $\mathbf{B}^q$ so that

$$(3.1) \qquad \phi_\nu(\zeta, v) = \Phi_n f_{\nu,n}(\zeta, v) = \sum_{\eta \in \mathbf{N}_0^q(n)} \psi_\nu^{(\eta)}(\zeta) \, \frac{\bar{v}^{\eta}}{\eta!},$$

as in (2.1). The functions $\psi_\nu^{(\eta)}$ were calculated in Corollary 2.6.

**3.1. Lemma.** *Let $\nu \in \mathbf{N}_0^q$ with $|\nu| \geq n$. Then there exists $\eta \in \mathbf{N}_0^q(n)$ such that $\psi_\nu^{(\eta)}$ is not identically zero on $\mathbf{B}^q$.*

*Proof.* Since $|\nu| \geq n$, there exists $\eta \in \mathbf{N}_0^q(n)$ such that $\nu \geq \eta$. By Corollary 2.6, $\psi_\nu^{(\eta)}$ is the desired nonzero component. ∎

**3.2. Lemma.** *Suppose that $\nu, \nu' \in \mathbf{N}_0^q$ with $|\nu'|, |\nu| \geq n$. Furthermore, suppose that there exists $\eta \in \mathbf{N}_0^q(n)$ with*

$$\psi_\nu^{(\eta)} = \psi_{\nu'}^{(\eta)} \neq 0.$$

*Then $\nu = \nu'$.*

*Proof.* By Corollary 2.6, $\psi_\nu^{(\eta)}(\zeta) = C_\nu \zeta^{\nu - \eta}$ and $\psi_{\nu'}^{(\eta)}(\zeta) = C_{\nu'} \zeta^{\nu' - \eta}$, where $C_\nu, C_{\nu'}$ are nonzero constants. Hence $\nu = \nu'$. ∎

Lemmas 3.1 and 3.2 show that $f_{\nu,n}$ can be recovered from a single nonzero summand $\psi_\nu^{(\eta)}(\zeta)\bar{v}^\eta/\eta!$ of $\phi_\nu(\zeta, v)$ in (3.1). However, there may be more than one nonzero summand of $\phi_\nu$. Our approach to the inversion problem is to first develop a method for isolating a nonzero summand of $\phi_\nu$.

We will construct a projection $P$ that when applied to $\phi_\nu$ returns the $\eta_\nu$-th summand $\psi_\nu^{(\eta_\nu)} \bar{v}^{\eta_\nu}/\eta_\nu!$ in (3.1), where $\eta_\nu \in \mathbf{N}_0^q(n)$ is maximal in the lexicographic order $\prec$ with respect to $\psi_\nu^{(\eta_\nu)}$ being nonzero. Restricting the $\omega_n$ action of Definition 1.10 to $K$, we see that the $K$ action on $\bar{\mathcal{P}}(n, \mathbf{C}^q)$ (thought of as the fiber over $\zeta = 0$ of the holomorphic, homogeneous vector bundle $\mathbf{E}_n$ for $G$ over $G/K$) is simply the natural action of $U(q)$ by left translation. Notice that each $\bar{v}^\eta$ is a weight vector for this representation, and that the lexicographic order $\prec$ corresponds to the order on these weight vectors, from lowest to highest. Hence $P$ will return the summand of $\phi_\nu$ of highest possible weight with respect to this action of $K$ on $\bar{\mathcal{P}}(n, \mathbf{C}^q)$.

**3.3. Definition.** For each $j \in \{1, 2, \ldots, q\}$ define $Q_j : \mathcal{O}(\mathbf{B}^q, \bar{\mathcal{P}}(n, \mathbf{C}^q)) \to \mathcal{O}(\mathbf{B}^q, \bar{\mathcal{P}}(n, \mathbf{C}^q))$ by

$$Q_j f(\zeta, v) = f(\zeta_1, \ldots, \zeta_{j-1}, 0, \zeta_{j+1}, \ldots, \zeta_q, v).$$

**3.4. Definition.** If $\eta \in \mathbf{N}_0^q(n)$ has an immediate successor $\tilde{\eta}$ with respect to the partial order $\prec$, let $s(\eta)$ be the index at which $\eta$ first disagrees with $\tilde{\eta}$. If $\eta$ has no immediate successor, i.e., $\eta = (n, 0, \ldots, 0)$, let $s(\eta) = 0$.

For example, $\eta = (3, 1, 2, 4, 0)$ has immediate successor $\tilde{\eta} = (3, 1, 3, 3, 0)$. So $s(\eta) = 3$, as $\eta$ and $\tilde{\eta}$ first disagree in the third component.

**3.5. Definition.** Suppose that $\eta \in \mathbf{N}_0^q(n)$. Define $P_\eta : \mathcal{O}(\mathbf{B}^q, \bar{\mathcal{P}}(n, \mathbf{C}^q)) \to \mathcal{O}(\mathbf{B}^q, \bar{\mathcal{P}}(n, \mathbf{C}^q))$ by

$$P_\eta := \prod_{j=1}^{s(\eta)} Q_j$$

where $s(\eta)$ is as in Definition 3.4. If $\eta = (n, 0, \ldots, 0)$, then $P_\eta := Id$.

**3.6. Definition.** Define $P : \mathcal{O}(\mathbf{B}^q, \bar{\mathcal{P}}(n, \mathbf{C}^q)) \to \mathcal{O}(\mathbf{B}^q, \bar{\mathcal{P}}(n, \mathbf{C}^q))$ by

$$(P\phi)(\zeta, v) = \sum_{\eta \in \mathbf{N}_0^q(n)} P_\eta \frac{\bar{v}^\eta}{\eta!} \frac{\partial^n \phi}{\partial \bar{v}^\eta}(\zeta, v).$$

The role of the operator $(\bar{v}^\eta/\eta!)(\partial^n/\partial\bar{v}^\eta)$ in Definition 3.6 is to single out the terms containing $\bar{v}^\eta$ in $\phi$, as in (2.1). Then $P_\eta$ is applied only to those terms.

For example, let $q = 3$ and $n = 2$. In this case the following table gives possible values of $\eta \in \mathbf{N}_0^3(2)$, in increasing lexicographic order, and the corresponding mappings $P_\eta$.

| $\eta$ | $s(\eta)$ | $P_\eta$ |
| --- | --- | --- |
| $(0,0,2)$ | 2 | $Q_1 Q_2$ |
| $(0,1,1)$ | 2 | $Q_1 Q_2$ |
| $(0,2,0)$ | 1 | $Q_1$ |
| $(1,0,1)$ | 2 | $Q_1 Q_2$ |
| $(1,1,0)$ | 1 | $Q_1$ |
| $(2,0,0)$ | 0 | $Id$ |

Any function $\phi \in \mathcal{O}(\mathbf{B}^3, \bar{\mathcal{P}}(2, \mathbf{C}^3))$ is written using (2.1) as

$$\phi(\zeta, v) = \psi^{(0,0,2)}(\zeta)\frac{\bar{v}_3^2}{2!} + \psi^{(0,1,1)}(\zeta)\bar{v}_2\bar{v}_3 + \cdots + \psi^{(2,0,0)}(\zeta)\frac{\bar{v}_1^2}{2!},$$

and the map $P$ acts by $P_\eta$, as given in the table, on the term containing $\bar{v}^\eta$.

Now let $f \in \mathcal{F}_2^{1,3}$ be given by $f(z) = z_1^3 \bar{z}_2 \bar{z}_3^3 \bar{z}_4$, so that $\nu = (1,3,1)$, in the notation of Corollary 2.6. The $\eta \in \mathbf{N}_0^3(2)$ for which $\nu \geq \eta$ are $(0,1,1)$, $(0,2,0)$, $(1,0,1)$, and $(1,1,0)$. It follows from Corollary 2.6 that

$$(\Phi_2 f)(\zeta, v) = 18\zeta_1\zeta_2^2 \bar{v}_2\bar{v}_3 + 36\zeta_1\zeta_2\zeta_3\frac{\bar{v}_2^2}{2!} + 6\zeta_2^3 \bar{v}_1 \bar{v}_3 + 18\zeta_2^2\zeta_3 \bar{v}_1 \bar{v}_2.$$

Using Definition 3.6 and the table above, we compute that for $\phi = \Phi_2 f$,

$$P\phi(\zeta, v) = 18\zeta_2^2\zeta_3 \bar{v}_1 \bar{v}_2.$$

Notice that only one term of $\phi(\zeta, v)$ remains, the term indexed by the maximal $\eta \in \mathbf{N}_0^3(2)$ with respect to lexicographic order for which $\psi^{(\eta)} \neq 0$. This is no coincidence, as is shown by the next lemma.

**3.7. Lemma.** *Fix $\nu \in \mathbf{N}_0^q$ with $|\nu| \geq n$. Let $f_{\nu,n}(z) = z_1^{|\nu|-n} \bar{z}_S^\nu$, and let $\Phi_n f_{\nu,n} = \phi_\nu = \sum_{\eta \in \mathbf{N}_0^q(n)} \psi_\nu^{(\eta)} \bar{v}^\eta/\eta!$, as in (3.1). Let $\alpha \in \mathbf{N}_0^q(n)$ be maximal in the order $\prec$ with respect to $\psi_\nu^{(\alpha)}$ being nonzero. Then $P_\alpha \psi_\nu^{(\alpha)} \bar{v}^\alpha/\alpha! = \psi_\nu^{(\alpha)} \bar{v}^\alpha/\alpha!$ and for $\eta \prec \alpha$, $P_\eta \psi_\nu^{(\eta)} \bar{v}^\eta/\eta! = 0$.*

*Proof.* Since $\psi_\nu^{(\eta)}$ is a monomial in $\zeta$, then either $P_\eta \psi_\nu^{(\eta)} \bar{v}^\eta/\eta! = \psi_\nu^{(\eta)} \bar{v}^\eta/\eta!$ or $P_\eta \psi_\nu^{(\eta)} \bar{v}^\eta/\eta! = 0$ for any $\eta \in \mathbf{N}_0^q(n)$. We first show that $P_\alpha \psi_\nu^{(\alpha)} \bar{v}^\alpha/\alpha! = \psi_\nu^{(\alpha)} \bar{v}^\alpha/\alpha!$.

Suppose that $\alpha = (n, 0, \ldots, 0)$. Then by Definition 3.5, $P_\alpha = Id$, so we have $P_\alpha \psi_\nu^{(\alpha)} \bar{v}^\alpha/\alpha! = \psi_\nu^{(\alpha)} \bar{v}^\alpha/\alpha!$ trivially.

So we now assume that $\alpha = (\alpha_1, \ldots, \alpha_q) \neq (n, 0, \ldots, 0)$. Let $\tilde{\alpha}$ be the immediate successor of $\alpha$ within $\mathbf{N}_0^q(n)$. Note that $\tilde{\alpha}_i = \alpha_i$ for $1 \leq i < s(\alpha)$. Also, $\tilde{\alpha}_{s(\alpha)} > \alpha_{s(\alpha)}$, and $|\tilde{\alpha}| = |\alpha| = n$. It follows that there must be a nonzero component $\alpha_t$ of $\alpha$ with $s(\alpha) < t \leq q$. Also, by Corollary 2.6, we have that

$$\psi_\nu^{(\alpha)}(\zeta) = C_\alpha \zeta_1^{\nu_1 - \alpha_1} \ldots \zeta_q^{\nu_q - \alpha_q},$$

where $\nu \geq \alpha$ and $C_\alpha$ is a constant. For a contradiction, assume that $P_\alpha \psi_\nu^{(\alpha)} \bar{v}^\alpha / \alpha! = 0$. Then, by Definition 3.5, there exists $k \in \{1, \ldots, s(\alpha)\}$ with $\nu_k - \alpha_k > 0$. Consider $\alpha' = \alpha - \epsilon_t + \epsilon_k$. We have that $\alpha' \in \mathbf{N}_0^q(n)$ with $\alpha' \succ \alpha$. Furthermore, by Corollary 2.6, $\psi_\nu^{(\alpha')} \neq 0$. This contradicts the maximality of $\alpha$.

Fix $\eta \prec \alpha$. The proof will be concluded by showing that $P_\eta \psi_\nu^{(\eta)} \bar{v}^\eta / \eta! = 0$. If $\psi_\nu^{(\eta)} = 0$ then the result is immediate. So, for a contradiction, assume that $\psi_\nu^{(\eta)} \neq 0$ and that $P_\eta \psi_\nu^{(\eta)} \bar{v}^\eta / \eta! = \psi_\nu^{(\eta)} \bar{v}^\eta / \eta!$. Definition 3.5 and Corollary 2.6 then imply that $\nu_i - \eta_i = 0$ for $i \in \{1, \ldots, s(\eta)\}$. Now, since $\alpha \succ \eta$, there exists $k \in \{1, \ldots, s(\eta)\}$ such that $\alpha_k > \eta_k$. However, then $\nu_k - \alpha_k < 0$, and so by Corollary 2.6, $\psi_\nu^{(\alpha)} = 0$. This contradicts the fact that $\psi_\nu^{(\alpha)} \neq 0$. ∎

We now construct another mapping of $\mathcal{O}(\mathbf{B}^q, \bar{\mathcal{P}}(n, \mathbf{C}^q))$ into itself that we will use to weight the mapping $P$. Essentially, we will be assigning to each element $\phi \in \mathcal{O}(\mathbf{B}^q, \bar{\mathcal{P}}(n, \mathbf{C}^q))$ a solution of a certain partial differential equation.

**3.8. Definition.** Suppose $f \in \mathcal{O}(\mathbf{B}^q, \bar{\mathcal{P}}(n, \mathbf{C}^q))$, $\zeta \in \mathbf{B}^q$, and $v \in \mathbf{C}^q$. For a given non-negative integer $m$ and an index $j \in \{1, \ldots, q\}$, we define the mapping $I(m, j) : \mathcal{O}(\mathbf{B}^q, \bar{\mathcal{P}}(n, \mathbf{C}^q)) \to \mathcal{O}(\mathbf{B}^q, \bar{\mathcal{P}}(n, \mathbf{C}^q))$ by

$$(I(m,j)f)(\zeta, v)$$
$$:= \int_0^1 \int_0^{t_m} \cdots \int_0^{t_2} f(\zeta_1, \ldots, \zeta_{j-1}, t_1 \zeta_j, \zeta_{j+1}, \ldots, \zeta_q, v) \, dt_1 dt_2 \ldots dt_m,$$

with the convention that $(I(0, j)f)(\zeta, v) = f(\zeta, v)$.

For example, if $f(\zeta, v) = \zeta_j^k \bar{v}_1^n$, then $(I(m, j)f)(\zeta, v) = \dfrac{k!}{(k+m)!} \zeta_j^k \bar{v}_1^n$.

**3.9. Definition.** Define $F_\eta : \mathcal{O}(\mathbf{B}^q, \bar{\mathcal{P}}(n, \mathbf{C}^q)) \to \mathcal{O}(\mathbf{B}^q, \bar{\mathcal{P}}(n, \mathbf{C}^q))$ by

$$(F_\eta f)(\zeta, v) := \Big(I(\eta_1, 1) I(\eta_2, 2) \ldots I(\eta_q, q) f\Big)(\zeta, v),$$

where $\eta \in \mathbf{N}_0^q(n)$.

Suppose that $f \in \mathcal{O}(\mathbf{B}^q, \bar{\mathcal{P}}(n, \mathbf{C}^q))$ and

$$f(\zeta, v) = \sum_{\eta \in \mathbf{N}_0^q(n)} \sum_{\gamma \in \mathbf{N}_0^q} a_{\eta, \gamma} \zeta^\gamma \bar{v}^\eta$$

is the Taylor series expansion for $f$ centered at the origin. Since the series above converges uniformly on compact subsets of $\mathbf{B}^q$, it follows from Definition 3.8 that if $\rho \in \mathbf{N}_0^q(n)$,

$$F_\rho f(\zeta, v) = \sum_{\eta \in \mathbf{N}_0^q(n)} \sum_{\gamma \in \mathbf{N}_0^q} a_{\eta,\gamma} F_\rho(\zeta^\gamma \bar{v}^\eta)$$

(3.2)
$$= \sum_{\eta \in \mathbf{N}_0^q(n)} \sum_{\gamma \in \mathbf{N}_0^q} \frac{a_{\eta,\gamma} \gamma!}{(\gamma + \rho)!} \zeta^\gamma \bar{v}^\eta.$$

Now, we weight the mapping $P$ by the $F_\eta$ to obtain a new mapping that will play a large role in the rest of the paper.

**3.10. Definition.** Define $L : \mathcal{O}(\mathbf{B}^q, \bar{\mathcal{P}}(n, \mathbf{C}^q)) \to \mathcal{O}(\mathbf{B}^q, \bar{\mathcal{P}}(n, \mathbf{C}^q))$ by

$$(L\phi)(\zeta, v) = \sum_{\eta \in \mathbf{N}_0^q(n)} F_\eta P_\eta \bar{v}^\eta \frac{\partial^n \phi}{\partial \bar{v}^\eta}(\zeta, v).$$

For example, recall that if $f(z) = z_1^3 \bar{z}_2 \bar{z}_3^3 \bar{z}_4$, and $\phi = \Phi_2 f$, then

$$P\phi(\zeta, v) = 18 \zeta_2^2 \zeta_3 \bar{v}_1 \bar{v}_2.$$

Applying (3.2) and Definition 3.10, we compute that

$$L\phi(\zeta, v) = 18 \frac{(0, 2, 1)!}{((0, 2, 1) + (1, 1, 0))!} \zeta_2^2 \zeta_3 \bar{v}_1 \bar{v}_2 = 6 \zeta_2^2 \zeta_3 \bar{v}_1 \bar{v}_2.$$

For $\phi \in \mathcal{O}(\mathbf{B}^q, \bar{\mathcal{P}}(n, \mathbf{C}^q))$, we express $\phi$ as the sum

$$\phi(\zeta, v) = \sum_{\eta \in \mathbf{N}_0^q(n)} \psi^{(\eta)}(\zeta) \frac{\bar{v}^\eta}{\eta!}$$

as in (2.1). For a compact subset $A$ of $\mathbf{B}^q$, set

(3.3)
$$\|\phi\|_A^2 := \sum_{\eta \in \mathbf{N}_0^q(n)} |\psi^{(\eta)}|_A^2,$$

where $|\cdot|_A$ is the supremum norm over $A$. We endow $\mathcal{O}(\mathbf{B}^q, \bar{\mathcal{P}}(n, \mathbf{C}^q))$ with the topology of uniform convergence on compact subsets of $\mathbf{B}^q$, i.e., the topology induced by the norm of (3.3).

**3.11. Lemma.** *The mapping $L : \mathcal{O}(\mathbf{B}^q, \bar{\mathcal{P}}(n, \mathbf{C}^q)) \to \mathcal{O}(\mathbf{B}^q, \bar{\mathcal{P}}(n, \mathbf{C}^q))$ is continuous.*

*Proof.* The remarks after Definition 3.9 serve to show that $F_\eta$ is continuous for each $\eta \in \mathbf{N}_0^q(n)$. It remains to show that $P_\eta$ is continuous.

Let $\phi \in \mathcal{O}(\mathbf{B}^q, \bar{\mathcal{P}}(n, \mathbf{C}^q))$. Suppose that $\{\phi_k\}$ is a sequence of functions in $\mathcal{O}(\mathbf{B}^q, \bar{\mathcal{P}}(n, \mathbf{C}^q))$ such that $\{\phi_k\}$ converges to $\phi$ uniformly on compact subsets of $\mathbf{B}^q$, in the norm of (3.3). Let $A$ be compact in $\mathbf{B}^q$ and for $1 \leq j \leq q$, let $A_j$ be a compact subset of $\mathbf{B}^q$ containing the natural projection of $A$ onto the set $\{\zeta \in \mathbf{B}^q \mid \zeta_j = 0\}$. Then

$$\|Q_j \phi - Q_j \phi_k\|_A \leq \|\phi - \phi_k\|_{A_j}.$$

It follows that $Q_j$ is continuous and hence $P_\eta$ is continuous. ∎

We are now ready to invert the transform $\Phi_n$ of section 2.

**3.12. Theorem.** *Suppose that $\phi \in \Phi_n(\mathcal{F}_n^{1,q})$, and that $\phi = \Phi_n f$, for $f \in \mathcal{F}_n^{1,q}$. Then, for $z \in \mathbf{C}^{q+1}$,*

$$f(z) = \lim_{t \to 1^-} \frac{d^q}{dt^q}\left[t^q \int_{\mathbf{B}^q}\int_{\mathbf{C}^q} (L\phi)(t\zeta, v)e^{(z_1\bar{\zeta}+v^T)\bar{z}_S}\, e^{-|v|^2}\, dm(v)dm(\zeta)\right].$$

*Proof.* The theorem is proven in two steps. First, given $\nu \in \mathbf{N}_0^q$ with $|\nu| \geq n$, we let $f_{\nu,n}(z) = z_1^{|\nu|-n}\bar{z}_S^\nu$ and we let $\phi_\nu = \Phi_n f_{\nu,n}$. We will prove that the theorem holds for $\phi_\nu$. Then we will show that the theorem holds for an arbitrary $\phi \in \Phi_n(\mathcal{F}_n^{1,q})$.

So, let $\nu$ be as above, and let $\eta_\nu \in \mathbf{N}_0^q(n)$ be such that $P_{\eta_\nu}(\psi_\nu^{(\eta_\nu)}\bar{v}^{\eta_\nu}/\eta_\nu!) \neq 0$. Then by Lemma 3.7 and Definition 3.10,

$$(L\phi_\nu)(\zeta, v) = \sum_{\eta \in \mathbf{N}_0^q(n)} F_\eta P_\eta \bar{v}^\eta \frac{\partial^n \phi_\nu}{\partial \bar{v}^\eta}(\zeta, v)$$

$$= F_{\eta_\nu} P_{\eta_\nu} \frac{\nu!(|\nu|-n)!}{(\nu-\eta_\nu)!} \zeta^{\nu-\eta_\nu}\bar{v}^{\eta_\nu}$$

(3.4) $$= (|\nu|-n)!\zeta^{\nu-\eta_\nu}\bar{v}^{\eta_\nu}.$$

Hence, by using (3.4), (1.3), and Lemma 1.9, we compute that

$$t^q \int_{\mathbf{B}^q}\int_{\mathbf{C}^q} (L\phi_\nu)(t\zeta, v)e^{(z_1\bar{\zeta}+v^T)\bar{z}_S}\, e^{-|v|^2}\, dm(v)dm(\zeta)$$

$$= t^{|\nu|-n+q}(|\nu|-n)! \int_{\mathbf{B}^q}\int_{\mathbf{C}^q} \zeta^{\nu-\eta_\nu}\bar{v}^{\eta_\nu} e^{(z_1\bar{\zeta}+v^T)\bar{z}_S}\, e^{-|v|^2}\, dm(v)dm(\zeta)$$

$$= t^{|\nu|-n+q}\bar{z}_S^{\eta_\nu}(|\nu|-n)! \int_{\mathbf{B}^q} \zeta^{\nu-\eta_\nu} \exp\left[z_1\left(\sum_{j=1}^q \overline{z_{j+1}\zeta_j}\right)\right] dm(\zeta)$$

$$= t^{|\nu|-n+q}\left[\frac{z_1^{|\nu|-n}\bar{z}_S^\nu(|\nu|-n)!}{(\nu-\eta_\nu)!}\right]\int_{\mathbf{B}^q} \zeta^{\nu-\eta_\nu}\bar{\zeta}^{\nu-\eta_\nu}\, dm(\zeta)$$

$$= t^{|\nu|-n+q}\left[\frac{z_1^{|\nu|-n}\bar{z}_S^\nu(|\nu|-n)!}{(|\nu|-n+q)!}\right],$$

and so

$$\lim_{t \to 1^-} \frac{d^q}{dt^q}\left[t^q \int_{\mathbf{B}^q}\int_{\mathbf{C}^q} (L\phi_\nu)(t\zeta, v)e^{(z_1\bar{\zeta}+v^T)\bar{z}_S}\, e^{-|v|^2}\, dm(v)dm(\zeta)\right]$$

$$= \lim_{t \to 1^-}\left(\frac{d^q}{dt^q}t^{|\nu|-n+q}\left[\frac{z_1^{|\nu|-n}\bar{z}_S^\nu(|\nu|-n)!}{(|\nu|-n+q)!}\right]\right)$$

$$= f_{\nu,n}(z).$$

This completes the first part of the proof.

Now let $\phi$ be an arbitrary element of $\Phi_n(\mathcal{F}_n^{1,q})$, $\phi = \Phi_n f$ for $f \in \mathcal{F}_n^{1,q}$. Put

$$f(z) = \sum_{\substack{\nu \in \mathbf{N}_0^q \\ |\nu| \geq n}} a_\nu f_{\nu,n}(z).$$

This series converges uniformly on compact subsets of $\mathbf{C}^{q+1}$ *and in the norm on* $\mathcal{F}_n^{1,q}$ *induced by* $\langle \cdot, \cdot \rangle$. Hence by Lemma 2.3, given $\zeta \in \mathbf{B}^q$,

$$\phi(\zeta, v) = \sum_{\substack{\nu \in \mathbf{N}_0^q \\ |\nu| \geq n}} a_\nu \phi_\nu(\zeta, v). \tag{3.5}$$

Moreover, the series in (3.5) constitutes a Taylor series expansion for $\phi$, and so the series converges uniformly on compact subsets of $\mathbf{B}^q$. Hence by (1.3) and Lemma 3.11, for $t \in [0,1)$ we have

$$\int_{\mathbf{B}^q} \int_{\mathbf{C}^q} L\bigg(\sum_{\substack{\nu \in \mathbf{N}_0^q \\ |\nu| \geq n}} a_\nu \phi_\nu\bigg)(t\zeta, v) e^{(z_1 \bar{\zeta} + v^T)\bar{z}_S} \, e^{-|v|^2} \, dm(v) dm(\zeta)$$

$$= \int_{\mathbf{B}^q} \int_{\mathbf{C}^q} \sum_{\substack{\nu \in \mathbf{N}_0^q \\ |\nu| \geq n}} a_\nu (L\phi_\nu)(t\zeta, v) e^{(z_1 \bar{\zeta} + v^T)\bar{z}_S} \, e^{-|v|^2} \, dm(v) dm(\zeta)$$

$$= \int_{\mathbf{B}^q} \sum_{\substack{\nu \in \mathbf{N}_0^q \\ |\nu| \geq n}} t^{|\nu|-n} a_\nu (|\nu|-n)! \zeta^{\nu-\eta_\nu} \bar{z}_S^{\eta_\nu} \exp\bigg[z_1 \bigg(\sum_{j=1}^q \overline{z_{j+1}\zeta_j}\bigg)\bigg] dm(\zeta). \tag{3.6}$$

Notice that we may exchange integration over $\mathbf{C}^q$ with summation in the equation preceding (3.6) due to the absolute convergence of the integral in $v$. We finish by showing that we may interchange integration and summation in the right hand side of (3.6). The theorem will then follow immediately from the first part of the proof. Thus by the Fubini-Tonelli theorem, it suffices to show that

$$\sum_{\substack{\nu \in \mathbf{N}_0^q \\ |\nu| \geq n}} |a_\nu|(|\nu|-n)! t^{|\nu|-n} \int_{\mathbf{B}^q} \bigg|\zeta^{\nu-\eta_\nu} \bar{z}_S^{\eta_\nu} \exp\bigg[z_1\bigg(\sum_{j=1}^q \overline{z_{j+1}\zeta_j}\bigg)\bigg]\bigg| dm(\zeta) < \infty.$$

Note that for a fixed $z \in \mathbf{C}^{q+1}$, the function

$$\zeta \mapsto \bar{z}_S^{\eta_\nu} \exp\bigg[z_1\bigg(\sum_{j=1}^q \overline{z_{j+1}\zeta_j}\bigg)\bigg]$$

is bounded in modulus on $\mathbf{B}^q$. Hence there is a positive constant $c_z$ such that

$$\sum_{\substack{\nu \in \mathbf{N}_0^q \\ |\nu| \geq n}} |a_\nu|(|\nu|-n)! t^{|\nu|-n} \int_{\mathbf{B}^q} \bigg|\zeta^{\nu-\eta_\nu} \bar{z}_S^{\eta_\nu} \exp\bigg[z_1\bigg(\sum_{j=1}^q \overline{z_{j+1}\zeta_j}\bigg)\bigg]\bigg| dm(\zeta)$$

$$\leq c_z \sum_{\substack{\nu \in \mathbf{N}_0^q \\ |\nu| \geq n}} |a_\nu|(|\nu|-n)! t^{|\nu|-n} \int_{\mathbf{B}^q} |\zeta^{\nu-\eta_\nu}| dm(\zeta). \tag{3.7}$$

Now by applying the Cauchy-Schwarz inequality to the integral in (3.7) and using Lemma 1.9 we compute that

$$c_z \sum_{\substack{\nu \in \mathbf{N}_0^q \\ |\nu| \geq n}} |a_\nu|(|\nu|-n)! t^{|\nu|-n} \int_{\mathbf{B}^q} |\zeta^{\nu-\eta_\nu}| \, dm(\zeta)$$

$$\leq c_z' \sum_{\substack{\nu \in \mathbf{N}_0^q \\ |\nu| \geq n}} |a_\nu|(|\nu|-n)! t^{|\nu|-n} \Big(\frac{(\nu-\eta_\nu)!}{(|\nu|-n+q)!}\Big)^{\frac{1}{2}}$$

$$\leq c_z' \sum_{\substack{\nu \in \mathbf{N}_0^q \\ |\nu| \geq n}} |a_\nu|((|\nu|-n)!)^{\frac{1}{2}}(\nu!)^{\frac{1}{2}}\Big(\frac{t^{2(|\nu|-n)}}{(|\nu|-n+q)\ldots(|\nu|-n+1)}\Big)^{\frac{1}{2}}$$

$$(3.8) \qquad \leq c_z' \Bigg[\sum_{\substack{\nu \in \mathbf{N}_0^q \\ |\nu| \geq n}} |a_\nu|^2 (|\nu|-n)! \nu!\Bigg]^{\frac{1}{2}} \Bigg[\sum_{\substack{\nu \in \mathbf{N}_0^q \\ |\nu| \geq n}} \frac{t^{2(|\nu|-n)}}{(|\nu|-n+q)\ldots(|\nu|-n+1)}\Bigg]^{\frac{1}{2}}.$$

By (1.2), the first factor of (3.8) is simply $\langle f, f \rangle^{\frac{1}{2}}$, while the second factor is finite by Lemma 1.4. Hence the product in (3.8) is finite, and so the proof is concluded. ∎

## 4. A Hilbert Space Structure for $\Phi_n(\mathcal{F}_n^{1,q})$

Recalling Theorems 1.7 and 2.2, observe that for $n \geq 0$ the space $\Phi_n(\mathcal{F}_n^{1,q})$ is an irreducible subspace of $\mathcal{O}(\mathbf{B}^q, \bar{\mathcal{P}}(n, \mathbf{C}^q))$ under the representation $\omega_n$ of $U(1,q)$. Let $\vartheta_n$ denote this subrepresentation. One would like to explicitly describe a Hilbert space structure on $\Phi_n(\mathcal{F}_n^{1,q})$ with the property that $\vartheta_n$ and $\sigma_n$ will be unitarily equivalent via $\Phi_n$. It is worth pointing out that $\sigma_n$ is not in the discrete series for any $n \geq 0$ if $p = 1$ and $q > 1$. If $p = q = 1$ then the ladder representation of spin $n \geq 1$ is in the discrete series and the representation with spin $n = 0$ is a limit of discrete series (see Theorem 4.4 of [3]).

**4.1. Definition.** For $\phi_1, \phi_2 \in \mathcal{O}(\mathbf{B}^q, \bar{\mathcal{P}}(n, \mathbf{C}^q))$, define

$$((\phi_1, \phi_2))_{n,q} := \lim_{t \to 1^-} \frac{d^q}{dt^q}\bigg[t^q \int_{\mathbf{B}^q} \int_{\mathbf{C}^q} (L\phi_1)(t\zeta, v)\overline{\phi_2(\zeta, v)}\, e^{-|v|^2}\, dm(v) dm(\zeta)\bigg]$$

$$= \lim_{t \to 1^-} \frac{d^q}{dt^q}\bigg[t^q \int_{\mathbf{B}^q} \langle (L\phi_1)(t\zeta, \cdot), \phi_2(\zeta, \cdot)\rangle\, dm(\zeta)\bigg],$$

whenever the integral above converges. Here $L$ is as in Definition 3.10.

For $\phi \in \mathcal{O}(\mathbf{B}^q, \bar{\mathcal{P}}(n, \mathbf{C}^q))$, put $\|\phi\|_{n,q}^2 := ((\phi, \phi))_{n,q}$. It is easy to show that on the set $\{\phi \in \mathcal{O}(\mathbf{B}^q, \bar{\mathcal{P}}(n, \mathbf{C}^q)) \mid \|\phi\|_{n,q}^2 < \infty\}$, the form $((\cdot, \cdot))_{n,q}$ is hermitian and positive semi-definite. We next examine its behavior on $\Phi_n(\mathcal{F}_n^{1,q})$.

**4.2. Theorem.** *Let $\phi_1, \phi_2 \in \Phi_n(\mathcal{F}_n^{1,q})$ with $\phi_1 = \Phi_n f_1$ and $\phi_2 = \Phi_n f_2$ for $f_1, f_2 \in \mathcal{F}_n^{1,q}$. Then*

$$((\phi_1, \phi_2))_{n,q} = \langle f_1, f_2 \rangle.$$

Hence $((\cdot, \cdot))_{n,q}$ defines a positive definite hermitian form on $\Phi_n(\mathcal{F}_n^{1,q})$.

*Proof.* Put

$$f_1(z) = \sum_{\substack{\nu \in \mathbf{N}_0^q \\ |\nu| \geq n}} a_\nu f_{\nu,n}(z) \quad \text{and} \quad f_2(z) = \sum_{\substack{\nu \in \mathbf{N}_0^q \\ |\nu| \geq n}} b_\nu f_{\nu,n}(z),$$

where $f_{\nu,n}$ is as in Corollary 2.6. Correspondingly, according to Lemma 2.3, we have

$$\phi_1(\zeta, v) = \sum_{\substack{\nu \in \mathbf{N}_0^q \\ |\nu| \geq n}} a_\nu \phi_\nu(\zeta, v) \quad \text{and} \quad \phi_2(\zeta, v) = \sum_{\substack{\nu \in \mathbf{N}_0^q \\ |\nu| \geq n}} b_\nu \phi_\nu(\zeta, v),$$

where $\phi_\nu = \Phi_n f_{\nu,n}$. Then by using the series expansions above, along with Lemma 3.11, Corollary 2.6, (3.4), and Lemma 1.2, for $t \in [0, 1)$,

$$\int_{\mathbf{B}^q} \int_{\mathbf{C}^q} (L\phi_1)(t\zeta, v)\overline{\phi_2(\zeta, v)} e^{-|v|^2} \, dm(v) dm(\zeta)$$

$$= \int_{\mathbf{B}^q} \int_{\mathbf{C}^q} \sum_{\substack{\nu \in \mathbf{N}_0^q \\ |\nu| \geq n}} \sum_{\substack{\tilde{\nu} \in \mathbf{N}_0^q \\ |\tilde{\nu}| \geq n}} a_\nu \bar{b}_{\tilde{\nu}} (L\phi_\nu)(t\zeta, v)\overline{\phi_{\tilde{\nu}}(\zeta, v)} \, e^{-|v|^2} \, dm(v) dm(\zeta)$$

$$(4.1) \qquad = \int_{\mathbf{B}^q} \left[ \sum_{\substack{\nu, \tilde{\nu} \in \mathbf{N}_0^q \\ |\nu|, |\tilde{\nu}| \geq n \\ \tilde{\nu} \geq \eta_\nu}} a_\nu \bar{b}_{\tilde{\nu}} t^{|\nu|-n} \frac{(|\nu|-n)!(|\tilde{\nu}|-n)!\tilde{\nu}!}{(\tilde{\nu}-\eta_\nu)!} \zeta^{\nu-\eta_\nu} \bar{\zeta}^{\tilde{\nu}-\eta_\nu} \right] dm(\zeta).$$

Since $t \in [0, 1)$, by the dominated convergence theorem we may interchange integration and summation in (4.1). After doing so and applying Lemma 1.9, the only nonzero terms remaining in the sum are those for which $\nu = \tilde{\nu}$, yielding that

$$\int_{\mathbf{B}^q} \left[ \sum_{\substack{\nu, \tilde{\nu} \in \mathbf{N}_0^q \\ |\nu|, |\tilde{\nu}| \geq n \\ \tilde{\nu} \geq \eta_\nu}} a_\nu \bar{b}_{\tilde{\nu}} t^{|\nu|-n} \frac{(|\nu|-n)!(|\tilde{\nu}|-n)!\tilde{\nu}!}{(\tilde{\nu}-\eta_\nu)!} \zeta^{\nu-\eta_\nu} \bar{\zeta}^{\tilde{\nu}-\eta_\nu} \right] dm(\zeta)$$

$$= \sum_{\substack{\nu \in \mathbf{N}_0^q \\ |\nu| \geq n}} a_\nu \bar{b}_\nu t^{|\nu|-n} [(|\nu|-n)!]^2 \frac{\nu!}{(|\nu|-n+q)!}$$

$$= \sum_{\substack{\nu \in \mathbf{N}_0^q \\ |\nu| \geq n}} \frac{a_\nu \bar{b}_\nu t^{|\nu|-n} \nu! (|\nu|-n)!}{(|\nu|-n+q)(|\nu|-n+(q-1))\ldots(|\nu|-n+1)}.$$

Hence

$$\lim_{t \to 1^-} \frac{d^q}{dt^q} \left[ t^q \int_{\mathbf{B}^q} \int_{\mathbf{C}^q} (L\phi_1)(t\zeta, v)\overline{\phi_2(\zeta, v)} \, e^{-|v|^2} \, dm(v) dm(\zeta) \right]$$

$$= \sum_{\substack{\nu \in \mathbf{N}_0^q \\ |\nu| \geq n}} a_\nu \bar{b}_\nu (|\nu|-n)! \nu!$$

$$= \langle f_1, f_2 \rangle. \quad \blacksquare$$

**4.3. Corollary.** *The vector space $\Phi_n(\mathcal{F}_n^{1,q})$, when endowed with the inner product $((\cdot, \cdot))_{n,q}$, is a Hilbert space.*

*Proof.* By Theorem 4.2, $\Phi_n(\mathcal{F}_n^{1,q})$ is a pre-Hilbert space. So all we need to show is that $\Phi_n(\mathcal{F}_n^{1,q})$ is complete. Suppose that $\{\phi_k\}$ is a sequence in $\Phi_n(\mathcal{F}_n^{1,q})$ that is Cauchy in the norm induced by $((\cdot, \cdot))_{n,q}$, and suppose that $\phi_k = \Phi_n(f_k)$ for $f_k \in \mathcal{F}_n^{1,q}$. Then by Theorem 4.2, the sequence $\{f_k\}$ is Cauchy in the Hilbert Space $\mathcal{F}_n^{1,q}$. Hence there is a function $f \in \mathcal{F}_n^{1,q}$ such that $f_k$ converges to $f$ in norm. Let $\phi = \Phi_n f$. Then $\phi \in \Phi_n(\mathcal{F}_n^{1,q})$ and $\phi_k$ converges to $\phi$ in norm. ∎

**4.4. Corollary.** *The representation $\vartheta_n$ is irreducible and unitary on the Hilbert space $\Phi_n(\mathcal{F}_n^{1,q})$ and is unitarily equivalent to $\sigma_n$ via $\Phi_n$.*

*Proof.* This follows from Theorems 2.2 and 4.2. ∎

At this point one may wonder whether the inner product $((\cdot, \cdot))_{n,q}$ characterizes $\Phi_n(\mathcal{F}_n^{1,q})$. In general, this inner product is not enough to characterize $\Phi_n(\mathcal{F}_n^{1,q})$. However, when we put the inner product together with the differential equations of Theorem 2.4, we are able to combine our results with those of [4] to obtain a complete characterization of $\Phi_n(\mathcal{F}_n^{1,q})$.

For $q > 1$, let $\mathcal{D}_{1,q}$ denote the set of differential operators consisting of all $2 \times 2$-minors of the matrix of differential operators

$$\square_{1,q} = \begin{pmatrix} \frac{\partial}{\partial \zeta_1} & \cdots & \frac{\partial}{\partial \zeta_q} \\ \frac{\partial}{\partial \bar{v}_1} & \cdots & \frac{\partial}{\partial \bar{v}_q} \end{pmatrix}$$

of Theorem 2.4. If $q = 1$, we will say that $\mathcal{D}_{1,1} = \emptyset$. Also, let

$$\mathcal{X}_n^{1,q} = \{\phi \in \mathcal{O}(\mathbf{B}^q, \bar{\mathcal{P}}(n, \mathbf{C}^q)) \mid \|\phi\|_{n,q}^2 < \infty \text{ and } D\phi = 0 \text{ for all } D \in \mathcal{D}_{1,q}\}.$$

**4.5. Theorem.** *Let $n \in \mathbf{N}_0$ and $q \in \mathbf{N}$. Then $\Phi_n(\mathcal{F}_n^{1,q}) = \mathcal{X}_n^{1,q}$.*

*Proof.* From Theorems 2.4 and 4.2, we know that $\Phi_n(\mathcal{F}_n^{1,q}) \subset \mathcal{X}_n^{1,q}$. So it remains to verify the reverse inclusion.

Suppose that $\phi \in \mathcal{X}_n^{1,q}$. Then

$$\phi(\zeta, v) = \sum_{m=0}^{\infty} \phi_{(m)}(\zeta, v),$$

where each $\phi_{(m)}$ is a $\bar{\mathcal{P}}(n, \mathbf{C}^q)$-valued holomorphic polynomial that is homogeneous of degree $m$ on $\mathbf{B}^q$. Since $0 = D\phi(\zeta, v) = \sum_{m=0}^{\infty} D\phi_{(m)}(\zeta, v)$, and since the operator $D$ preserves the homogeneity, hence the linear independence of the functions $D\phi_{(m)}$, we conclude that $D\phi_{(m)} = 0$ for all $m \in \mathbf{N}_0$. The proof of Theorem 4.2 implies that any polynomial in $\zeta$ has finite norm. Hence $\phi_{(m)} \in \mathcal{X}_n^{1,q}$ for each $m$.

In the case that $d(n, q) > 1$, the results of [4] imply that each $\phi_{(m)}$ is also an element of $\Phi_n(\mathcal{F}_n^{1,q})$. In fact, this also holds in the line bundle case, i.e., when $d(n, q) = 1$. Corollary 2.6 shows that when $d(n, q) = 1$, the image of a monomial is a monomial. Thus in this case we can invert $\Phi_n$ directly from Corollary 2.6.

So, with $\phi_{(m)}$ as above, there exists $f_{(m)} \in \mathcal{F}_n^{1,q}$ with $\Phi_n f_{(m)} = \phi_{(m)}$ for each $m \in \mathbf{N}_0$. By Corollary 2.6, $f_{(m)}$ must have the form

$$f_{(m)}(z) = \sum_{\substack{\nu \in \mathbf{N}_0^q \\ |\nu|=m+n}} a_{m,\nu} z_1^{|\nu|-n} \bar{z}_S^\nu.$$

Put $f = \sum f_{(m)}$. Since $\sum_{m=0}^N \phi_{(m)}$ converges to $\phi$ uniformly on compact subsets of $\mathbf{B}^q$ (see (3.3)), we may apply the dominated convergence theorem as in 4.2, along with Corollary 4.3, to get

$$\begin{aligned}
\|\phi\|_{n,q}^2 &= \lim_{N \to \infty} \|\sum_{m=0}^N \phi_{(m)}\|_{n,q}^2 \\
&= \lim_{N \to \infty} \langle \sum_{m=0}^N f_{(m)}, \sum_{m=0}^N f_{(m)} \rangle \\
&= \lim_{N \to \infty} \sum_{m=0}^N \sum_{\substack{\nu \in \mathbf{N}_0^q \\ |\nu|=m+n}} |a_{m,\nu}|^2 (|\nu|-n)! \nu! \\
&= \langle f, f \rangle.
\end{aligned}$$

Hence $\langle f, f \rangle < \infty$, and so $f \in \mathcal{F}_n^{1,q}$. From the continuity of $\Phi_n$ we then conclude that $\Phi_n f = \phi$.  ∎

Oklahoma State University, Stillwater, Oklahoma 74078–0613
*E-mail address*: lorch@math.okstate.edu

Oklahoma State University, Stillwater, Oklahoma 74078–0613
*Current address*: Institute for Advanced Study, Princeton, New Jersey 08540
*E-mail address*: mantini@math.okstate.edu